\documentclass[12pt]{article}
\usepackage{verbatim,amsmath,amssymb,amscd,amsthm,color,graphics,fancybox}
\usepackage[stable]{footmisc}
\usepackage[all]{xy}
\usepackage[margin=1.1in]{geometry}
\usepackage[enableskew]{youngtab}

\theoremstyle{definition}

\newcommand{\hako}{ 
\thicklines
\put(0,0){\line(0,1){0.8}}
\put(0,0){\line(1,0){0.8}}
\put(0.8,0){\line(0,1){0.8}}
\thinlines
}

\newcommand{\tama}{ 
\color[cmyk]{0,0,0,0.2}
\put(0.20,0.4){\circle*{0.6}}
\color{black}
\put(0.20,0.4){\circle{0.6}}
}

\newcommand{\kago}{ 
\put(0.05,0.3){\line(1,0){0.7}}
\qbezier(0.05,0.3)(0.05,0)(0.4,0)
\qbezier(0.4,0)(0.75,0)(0.75,0.3)
\qbezier(0.2,0.3)(0.2,0.2)(0.6,0.05)
\qbezier(0.4,0.3)(0.4,0.2)(0.7,0.15)
\qbezier(0.6,0.3)(0.6,0.2)(0.2,0.05)
\qbezier(0.4,0.3)(0.4,0.2)(0.1,0.15)
}

\begin{document}
\begin{center}
{\Huge \textsf{Ultradiscrete Soliton Systems and Combinatorial Representation Theory}}\\
\vspace{7mm}

{\Large \textsf{Reiho Sakamoto}}\\

\textsf{Department of Physics, Tokyo University of Science,\\ Kagurazaka Shinjukuku Tokyo Japan}
\end{center}

\section{Introduction}
This lecture note is intended to be a brief introduction to a recent development
on the interplay between the ultradiscrete (or tropical) soliton systems and
the combinatorial representation theory.
We will concentrate on the simplest cases which admit elementary
explanations without losing essential ideas of the theory.
In particular we give definitions for the main constructions corresponding
to the vector representation of type $A^{(1)}_1$.

This note is organized as follows.
In Section \ref{sec:bbs} we give a definition of the simplest example of the box-ball systems.
In Section \ref{sec:crystal} we explain a relationship between the
box-ball systems and the crystal bases of the quantum affine algebras.
In Section \ref{sec:RC} we give the definition of the rigged configuration bijection
for the vector representation of type $A^{(1)}_1$.
In Section \ref{sec:IST} we see that the rigged configurations give
the complete set of the action and angle variables for the box-ball systems.
This is the fundamental observation in the recent development on a
relationship between the box-ball systems and the combinatorial representation theory.
In Section \ref{sec:bbbs} we explain basic properties of the box-basket-ball systems
which are recently found generalizations of the box-ball systems.
The characteristic property of the system is that it is a mixture of
fermions and bosons with mutual interactions.
Finally, in Section \ref{sec:further}, we give comments on generalizations and further developments
of the materials discussed in this note.

This is the lecture note prepared for the conference ``Algebraic combinatorics related to
Young diagrams and statistical physics" held at International Institute for Advanced Studies (Kyoto)
during August 6--10, 2012.
The author is grateful to Professor Masao Ishikawa, Professor Soichi Okada
and Professor Hiroyuki Tagawa for the kind invitation to the conference and the warm hospitality.

\section{The box-ball systems}\label{sec:bbs}
In this section, let us define the simplest case of the {\bf box-ball systems}
introduced by Takahashi and Satsuma \cite{TS}.
The box-ball systems are prototypical examples of the ultradiscrete soliton systems.
Originally the ultradiscrete soliton system is a class of discrete dynamical systems
obtained by the ultradiscrete (or tropical) limit of the ordinary soliton systems \cite{TTMS}.
In this article we are interested in ultradiscrete soliton systems
which admit combinatorial interpretations.

Following the box and ball interpretation of the system \cite{Tak},
we prepare boxes which can accommodate at most one ball within each box.
We put many such boxes on a line and put finitely many balls of the same kind to the boxes.
We regard this configuration as the initial state of the system.
Then we perform the time evolution of the state by the following algorithm.
\begin{quotation}
\noindent
{\bf The time evolution of the box-ball system:}
Consider each ball from left to right and move the ball to the next available empty box.
Each ball is moved exactly once.
\end{quotation}
Here if necessary we put enough many empty boxes on the right of the given state
in order to keep the balls within the state.
We give an example of such time evolution starting from the top row
and proceeding downwards.
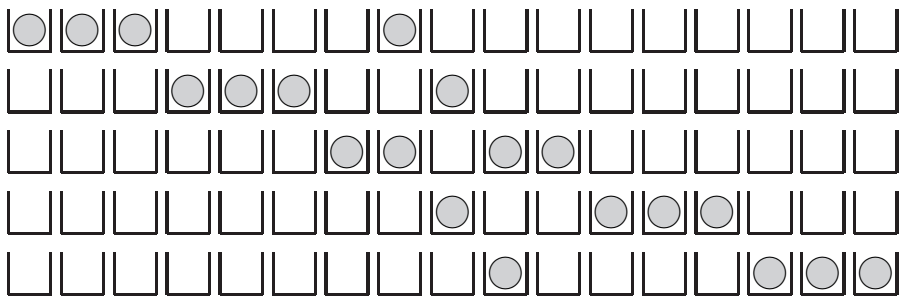
\begin{center}
\unitlength 20pt
\begin{picture}(17,1.1)
\put(0,0){\tama}
\put(1,0){\tama}
\put(2,0){\tama}
\put(7,0){\tama}
\multiput(0,0)(1,0){17}{\hako}
\end{picture}
\begin{picture}(17,1.1)
\put(3,0){\tama}
\put(4,0){\tama}
\put(5,0){\tama}
\put(8,0){\tama}
\multiput(0,0)(1,0){17}{\hako}
\end{picture}
\begin{picture}(17,1.1)
\put(6,0){\tama}
\put(7,0){\tama}
\put(9,0){\tama}
\put(10,0){\tama}
\multiput(0,0)(1,0){17}{\hako}
\end{picture}
\begin{picture}(17,1.1)
\put(8,0){\tama}
\put(11,0){\tama}
\put(12,0){\tama}
\put(13,0){\tama}
\multiput(0,0)(1,0){17}{\hako}
\end{picture}
\begin{picture}(17,1.1)
\put(9,0){\tama}
\put(14,0){\tama}
\put(15,0){\tama}
\put(16,0){\tama}
\multiput(0,0)(1,0){17}{\hako}
\end{picture}
\end{center}
In the box-ball system, we regard consecutive balls as solitary waves.
For example, the initial state in the above example contains two solitary waves
of length 3 and 1, respectively.
Then our interpretation of the above example is as follows.
If there is no interaction between waves, they move at velocity equal to each length.
In the course of the time evolution, solitary waves make collision with each other,
though they retain their original shapes after the collision except for the changes in the
positions compared with the possible positions we would have
if there is no interaction between waves.
Such properties of the waves of the box-ball systems
are characteristic of the soliton systems (see, e.g., \cite{kobo})
and we will call such waves {\bf solitons}.

Here we give a short list of remarks on the early papers.
During 1980's, there were several attempts of finding cellular automata
with solitonic properties.
A typical example of such researches is the filter automata
introduced by Park, Steiglitz and Thurston \cite{PST}.
In 1990, Takahashi and Satsuma \cite{TS} introduced the simplest case of the
box-ball systems and Takahashi \cite{Tak} described the algorithm
in terms of the boxes and balls.
The above definition corresponds to the original Takahashi--Satsuma box ball system.
A relationship between the Takahashi--Satsuma box-ball system and
the ordinary soliton systems including the KdV equation is discovered by Tokihiro, Takahashi,
Matsukidaira and Satsuma \cite{TTMS} via the limiting procedure
called the ultradiscrete limit.
A connection with the Toda equation is discussed in \cite{NTT}.
Such connections between the box-ball system and the usual soliton systems
show the classical integrability of the box-ball system.

Takahashi's box and ball algorithm provides several generalizations of the original box-ball system.
For example, in \cite{Tak} an internal degree of freedom for the balls (balls with different colors)
is introduced.
A connection with the Toda equation in the generalized context is discussed in \cite{TNS}.
The other degrees of freedom called a carrier \cite{TM} or a capacity of boxes \cite{TTM} are also introduced.
Such combinatorial interpretations of the time evolutions give nice intuition about the models
in many cases.

\section{A connection with the crystal bases}\label{sec:crystal}
A very important fact \cite{HHIKTT,FOY} about the box-ball systems is that their dynamics is
in fact governed by the Kashiwara's {\bf crystal bases} \cite{Kas} for the quantum affine algebras.
This formalism includes all extensions of the box-ball system which are mentioned in the last section.
Although the formulation does not depend on the types of the algebra,
we will concentrate on the simplest case $A^{(1)}_1$ here.

In order to describe the formulation, we need to consider more general boxes
which have capacities more than one.
Let $(a,b)$ represents the box of capacity $a+b$ containing $b$ balls.
Then the state $(a,b)$ can accommodate extra $a$ balls.
Let us denote the set of all such states as
\begin{align}
B^{1,s}:=\bigl\{(a,b)\,|\,a,b\in\mathbb{Z}_{\geq 0},\,a+b=s\bigr\}
\end{align}
which we call {\bf crystals}.\footnote{In general, we can identify the Kirillov--Reshetikhin
crystals $B^{r,s}$ for type $A^{(1)}_{n-1}$ with the set of $r\times s$ semistandard tableaux
with letters $1,2,\ldots,n$.
In this identification our $(a,b)$ is the height one semistandard tableau
with $a$ 1's and $b$ 2's. $B^{r,s}$ corresponds to the Kirillov--Reshetikhin module
naturally corresponding to the weight $s\Lambda_r$ where $\Lambda_r$ is the $r$-th fundamental weight.}
In particular we call $B^{1,1}$ the crystals for the vector representation.
In this coordinate, the states of the box-ball system in the previous section
are sequences of balls $(0,1)$ and empty places $(1,0)$.
Then we represent the states as $(1,0)\otimes (0,1)\otimes\cdots$
where $\otimes$ is the tensor product of crystals
(the readers may regard this as just alternative notation).
We call such elements of tensor products {\bf paths}.

The main ingredient of the formalism is the map called the {\bf combinatorial $R$-matrices}
\begin{align}
\begin{array}{llll}
R:&B^{1,s}\otimes B^{1,s'}&\longrightarrow &B^{1,s'}\otimes B^{1,s}\\
   &(a,b)\otimes (c,d)&\longmapsto            &(c',d')\otimes (a',b').
\end{array}
\end{align}
In the present case $A^{(1)}_1$, the explicit form of the map is
\begin{align}
a'&= a + \min(b, c) - \min(a, d)\nonumber\\
b'&= b - \min(b, c) + \min(a, d)\nonumber\\
c'&= c - \min(b, c) + \min(a, d)\nonumber\\
d'&= d + \min(b, c) - \min(a, d).
\end{align}
An important point of the map $R$ is that it has a deep mathematical origin
as the intertwining map that interchanges left and right of the tensor products of crystals.
For the later purposes we introduce a vertex diagram
for the map $R:a\otimes b\mapsto b'\otimes a'$ as follows:
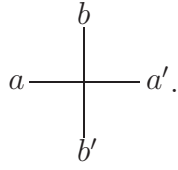
\begin{center}
\unitlength 13pt
\begin{picture}(4,4.5)
\put(0,2.0){\line(1,0){3.2}}
\put(1.6,0.4){\line(0,1){3.2}}
\put(-0.6,1.8){$a$}
\put(1.4,-0.3){$b'$}
\put(1.4,3.7){$b$}
\put(3.4,1.8){$a'$}
\put(4.1,1.7){.}
\end{picture}
\end{center}

By a repeated use of the map $R$ we define the time evolution of the box-ball systems
$T_l$ $(l\in{\mathbb Z}_{\geq 1})$ as follows.
Let $u_l:=(l,0)$ be the empty box of capacity $l$ and
let $b=b_1\otimes b_2\otimes\cdots\otimes b_L$ be a given state of the box-ball system.
We call $u_l$ the {\bf carrier}.
If necessary we put enough many empty boxes $(1,0)$ on the right.
Then we define $b_1'$, $\ldots$, $b_L'$ by the following diagram.
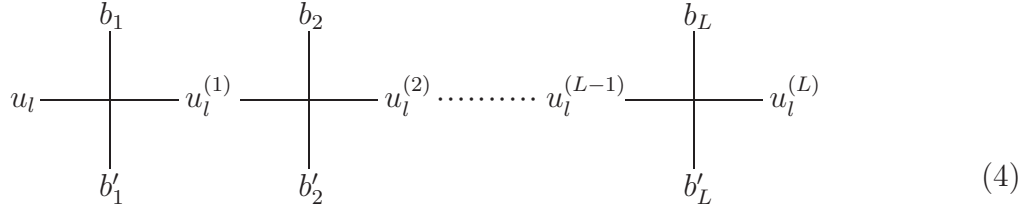
\begin{equation}\label{def:T_l}
\unitlength 13pt
\begin{picture}(22,5)(0,-0.5)
\multiput(0,0)(5.8,0){2}{
\put(0,2.0){\line(1,0){4}}
\put(2,0){\line(0,1){4}}
}
\put(-0.9,1.8){$u_l$}
\put(1.7,4.2){$b_1$}
\put(1.7,-0.8){$b_1'$}
\put(4.2,1.8){$u_l^{(1)}$}
\put(7.5,4.2){$b_2$}
\put(7.5,-0.8){$b_2'$}
\put(10.0,1.8){$u_l^{(2)}$}
\multiput(11.5,1.8)(0.3,0){10}{$\cdot$}
\put(14.7,1.8){$u_l^{(L-1)}$}
\put(17,0){
\put(0,2.0){\line(1,0){4}}
\put(2,0){\line(0,1){4}}
}
\put(18.7,4.2){$b_L$}
\put(18.7,-0.8){$b_L'$}
\put(21.2,1.8){$u_l^{(L)}$}
\end{picture}
\end{equation}
Here the precise meaning of the diagram is as follows.
We compute $R:u_l\otimes b_1\mapsto b_1'\otimes u_l^{(1)}$.
Then by using $u_l^{(1)}$ we compute
$R:u_l^{(1)}\otimes b_2\mapsto b_2'\otimes u_l^{(2)}$.
We do this procedure recursively until the end of the state.
Then we define
\begin{align}
T_l(b):=b_1'\otimes b_2'\otimes\cdots\otimes b_L'.
\end{align}
We can see that the time evolution rule given in Section \ref{sec:bbs}
coincides with $T_\infty$ here.

As a benefit of the definition by the crystal bases,
we can show the quantum integrability of the box-ball system
as the consequence of the Yang--Baxter relation for the combinatorial $R$-matrices \cite{FOY}.
More precisely, we have
\begin{align}
T_lT_k(b)=T_kT_l(b)
\end{align}
for arbitrary $l,k\in\mathbb{Z}_{\geq 1}$
and states $b$.
Moreover, we can construct conserved quantities of the box-ball system as follows.
Let us define (see(\ref{def:T_l}))
\begin{align}\label{def:E_l}
E_l(b)&:=\sum_{i=1}^LH(u_l^{(i-1)}\otimes b_i),\qquad E_0(b):=0
\end{align}
where $u_l^{(0)}:=u_l$ and the {\bf energy function}
$H:B^{1,s}\otimes B^{1,s'}\longrightarrow\mathbb{Z}$ is defined by
\begin{align}
H\bigl((a,b)\otimes (c,d)\bigr):=\min(a,d).
\end{align}
Again an important point of the energy function is that it has
a deep mathematical origin and is the consequence of the infinite dimensional
symmetry of the quantum affine algebras.
Let us consider the affinization of the crystal $B$
\begin{equation}
\mathrm{Aff}(B)=\{b[d]\, |\, b\in B,\, d\in\mathbb{Z}\}.
\end{equation}
For elements of tensor products of $\mathrm{Aff}(B)$,
we introduce the {\bf affine combinatorial $R$-matrices} by
\begin{equation}\label{eq:affine}
R_{\mathit{aff}}:b_1[d_1]\otimes b_2[d_2]\longmapsto
b_2'[d_2-H(b_1\otimes b_2)]\otimes
b_1'[d_1+H(b_1\otimes b_2)],
\end{equation}
where we have $R:b_1\otimes b_2\mapsto b_2'\otimes b_1'$
under the combinatorial $R$ matrix.
Then by the Yang--Baxter relation for the affine combinatorial $R$-matrices we
see that $E_l$ are the conserved quantities of the box-ball systems \cite{FOY}:
\begin{align}
E_l(T_k(b))=E_l(b).
\end{align}

\section{The rigged configurations}\label{sec:RC}
Another important aspect of the box-ball systems is a connection with the rigged configurations.
In this section we give the definition of a special case of the rigged configuration bijection
corresponding to the vector representation of type $A^{(1)}_1$.
Although this case is simpler than the general case,
it is still nontrivial and we can see basic ideas of the theory.

Originally the rigged configurations are discovered through insightful analysis of
the Bethe ansatz for quantum integrable systems \cite{KKR,KR}.
The main ingredient of the theory is a bijection between the set of
rigged configurations and elements of the tensor products of crystals.
Such a bijection is generalized for highest weight elements of
tensor products of the arbitrary Kirillov--Reshetikhin crystals of type $A^{(1)}_n$
and its mathematical theory is established by an important
paper of Kirillov, Schilling and Shimozono \cite{KSS}.

In our case, a rigged configuration is composed of a Young diagram (called the {\bf configuration})
and integers (called the {\bf riggings}) associated with each row of the Young diagram.
Let $\nu_i$ ($i=1,\ldots,g$) be the lengths of the rows of the configuration
and let $J_i$ be the rigging associated with the row $\nu_i$.
Then we represent the rigged configuration as $(\nu,J)=\{(\nu_i,J_i)\}_{i=1}^g$.
We call each $(\nu_i,J_i)$ {\bf string}.
Although there is a characterization of the possible rigged configurations,
we regard the set of the rigged configurations as the set of objects
obtained by the map (in fact, bijection)
\begin{align}
\Phi:b\longmapsto (\nu,J)
\end{align}
from arbitrary paths $b$.
We call the bijection $\Phi$ the {\bf rigged configuration bijection}.

Let us define the algorithm of the bijection $\Phi$.
For the given Young diagram $\nu$ let $Q_\ell(\nu)$
be the number of boxes contained in the left $\ell$ columns of $\nu$.
Suppose that we are given the path $b=b_1\otimes b_2\otimes\cdots\otimes b_L\in (B^{1,1})^{\otimes L}$
where the positions of the balls $b_k=(0,1)$
are given by $k=k_1,k_2,\ldots$ from left to right.
Let $P_\ell(k,\nu)$ be the {\bf vacancy number} defined by
\begin{equation}\label{def:P}
P_\ell(k,\nu):=k-2Q_\ell(\nu).
\end{equation}
For example, we have $P_2(16,\,\Yvcentermath1\Yboxdim5pt\yng(4,3,1)\,)=16-2\cdot 5=6$.
Suppose that a length $L$ path $b$ corresponds to the rigged configuration $(\nu,J)$.
Then we call the string $(\nu_i,J_i)$ {\bf singular} if the rigging $J_i$ coincides with
the corresponding vacancy number, that is, $P_{\nu_i}(L,\nu)=J_i$.

The bijection $\Phi$ is defined by a recursive procedure
corresponding to the positions of balls $k_1,k_2,\ldots$.
We start from the empty rigged configuration.
\begin{enumerate}
\item
Suppose that we have done the procedure up to $k_{j-1}$
and obtained the intermediate rigged configuration $(\eta,I)$.
\item
For the next position $k_j$, we do the following.
Suppose that the rigged configuration $(\eta,I)$ corresponds to a length $k_j-1$ path.
Compute the vacancy numbers $P_{\eta_i}(k_j-1,\eta)$ for all rows of $\eta$
and determine all the singular strings.
\item
If there is no singular string, add a length one row to the bottom of $\eta$.
Otherwise choose one of the longest singular string and add a box to the corresponding row.
Denote by $\eta'$ the new configuration thus obtained.\footnote{The order of rows is not essential
in the definition of $\Phi$.}
\item
Define the new rigging $I'$ as follows.
For the strings that are not changed under $\eta\rightarrow\eta'$,
we choose the same riggings as before.
Let $\eta_i'$ be the changed row under $\eta\rightarrow\eta'$.
Then define the new rigging by $I_i'=P_{\eta_i'}(k_j,\eta')$ so that the string
$(\eta_i',I_i')$ is singular in $(\eta',I')$.
The output $(\eta',I')$ is the new rigged configuration
corresponding to the length $k_j$ path.
\item
Repeat the same procedure for all $k_j$.
Let $(\nu,J)$ be the final output.
Then define $\Phi(b)=(\nu,J)$.
\end{enumerate}
A Mathematica package for the above procedure is available at \cite{demo}.
If we reverse all the procedure we obtain the algorithm for $\Phi^{-1}$.
As examples, let us look at the example of the time evolution of
the box-ball system at Section \ref{sec:bbs}.
In the first line, the positions of balls $k_j$ are $1,2,3,8$.
Then the computation of $\Phi$ proceeds as follows:
\begin{align}
\emptyset\xrightarrow{\,\,1\,\,}\,\,
\Yvcentermath1\Yboxdim12pt
\yng(1)\,\,{-1}\xrightarrow{\,\,2\,\,}\,
\yng(2)\,\,{-2}\xrightarrow{\,\,3\,\,}\,
\yng(3)\,\,{-3}
\label{rei1}
\xrightarrow{\,\,8\,\,}\,
\unitlength 12pt
\begin{picture}(3,1)
\put(0,-0.7){\Yboxdim12pt\yng(3,1)}
\put(3.2,0.5){$-3$}
\put(1.2,-0.7){$4$}
\end{picture}
\end{align}
Here we put riggings on the right of the corresponding row and
put $k_j$ above the corresponding arrows.
Similarly, for the third line of the same example, we have
\begin{align}
\label{rei2}
\emptyset\xrightarrow{\,\,7\,\,}&\,\,
\Yvcentermath1\Yboxdim12pt
\yng(1)\,\,5\xrightarrow{\,\,8\,\,}\,
\yng(2)\,\,4\xrightarrow{\,\,10\,\,}\,
\unitlength 12pt
\begin{picture}(2.5,1)
\put(0,-0.7){\Yboxdim12pt\yng(2,1)}
\put(2.2,0.5){$4$}
\put(1.2,-0.7){$6$}
\end{picture}
\xrightarrow{\,\,11\,\,}\,
\unitlength 12pt
\begin{picture}(3.5,1)
\put(0,-0.7){\Yboxdim12pt\yng(3,1)}
\put(3.2,0.5){$3$}
\put(1.2,-0.7){$6$}
\end{picture}
\end{align}

\section{The inverse scattering formalism}\label{sec:IST}
The main observation on the relationship between the rigged configurations
and the box-ball systems is that the rigged configuration bijection gives
the inverse scattering formalism for the box-ball systems.
In order to get the ideas of the result, let us compare the two examples in (\ref{rei1}) and (\ref{rei2}).
Then we see that the shapes of the configurations are same and
the differences of the riggings are two times the lengths of the corresponding rows.
Here we have the factor 2 in the change of riggings since we apply $T_\infty$ twice.

In general, let $b$ be the given state and let $\Phi(b)=\{(\nu_i,J_i)\}_{i=1}^g$.
Then we have \cite{KOSTY}
\begin{equation}\label{eq:ist}
\Phi(T_l(b))=\{(\nu_i,J_i+\min(l,\nu_i))\}_{i=1}^g.
\end{equation}
This property is valid for general box-ball systems including
all cases that appeared in \cite{HHIKTT,FOY}.
The proof of this fact heavily relies on a deep theorem of
Kirillov--Schilling--Shimozono \cite{KSS}.\footnote{If two tensor products
$b$ and $b'$ are isomorphic under the combinatorial $R$-matrices $R:b\mapsto b'$,
we have $\Phi(b)=\Phi(b')$ \cite[Lemma 8.5]{KSS}. The proof depends on
a large part of the paper.}
Indeed, if we compare (\ref{rei1}) and (\ref{rei2})
we can see that this property is already nontrivial.
To summarize, configurations are the conserved quantities (action variables)
and the riggings are the linearlization parameter (angle variables) of the box-ball systems.
Since $\Phi$ is bijective, the rigged configurations give the complete set of the action
and angle variables of the box-ball systems.\footnote{In fact
if we restrict to consider the box-ball systems corresponding to 
the vector representation of type $A^{(1)}_1$,
we do not need to use heavy apparatus like rigged configurations.
For example, \cite{TTS} introduced a combinatorial method to obtain the conserved quantities.
In \cite{MIT}, a method to obtain the action and angle variables is derived, which
is shown to be the special case of the rigged configurations \cite{KS}.
In \cite{KOTY} it is conjectured that the rigged configurations give
the action and angle variables of the box-ball system corresponding to
the vector representation of type $A^{(1)}_1$.
This problem is considered in \cite{Takagi1} with differently defined bijection.
We remark that in \cite{F} the Robinson--Schensted--Knuth algorithm is used
to give some of the conserved quantities of the box-ball system corresponding
to the vector representations of type $A^{(1)}_n$
(so called $P$-symbols are conserved under the time evolutions).}

Once we know that the rigged configurations are the underlying mathematical structure
of the box-ball systems, we can prove several fundamental properties of the box-ball systems.
For example, the box-ball systems considered in \cite{HHIKTT,FOY} are shown to
be solitonic by introducing a method to explicitly extract solitons from paths
as elements of the affinization of the crystals \cite{Sak1}.
The main point of the proof of the result is to introduce a structure of
the affine combinatorial $R$-matrices on the rigged configurations
via careful combinatorial arguments.
We remark that the proof of the solitonic properties of the box-ball systems
corresponding to the vector representation of type $A^{(1)}_n$ is proved in \cite{TNS}
by taking certain ultradiscrete limit of an ordinary soliton system
and an elegant alternative proof of their result is given in \cite{FOY}
by using the crystal bases.

Another important problem that is solved by the rigged configuration bijection
is the initial value problem of the box-ball systems \cite{KSY}.
The result includes all the extensions considered in \cite{HHIKTT,FOY}.
We note that an equivalent result for the case of the vector representation of $A^{(1)}_1$
is rederived in \cite{MIT2}.
Let us explain the result for the case of the vector representation of $A^{(1)}_1$.
The main point is to give an explicit piecewise linear formula for the
map $\Phi^{-1}:(\nu,J)\longmapsto b$.
For the given rigged configuration $(\nu,J)=\{(\nu_i,J_i)\}_{i=1}^g$,
let us define the following {\bf ultradiscrete tau functions}:
\begin{align}\label{def:tau}
\tau_r(k) := &-\min_{n \,\in \{0,1\}^g}
\biggl\{\sum_{i=1}^g(J_i+r\nu_i -k)n_i 
+ \sum_{i,j =1}^g\min(\nu_i,\nu_j)n_in_j\biggr\},
\quad (r=0,1)
\end{align}
where we denote $n=(n_1,\ldots,n_g)$.\footnote{If we consider paths with
periodicities, these functions $\tau_r$ exactly coincide with the {\bf tropical Riemann theta function} \cite{KuS}.}
Let us represent the $k$-th element of the path $b$ as $b_k=(1-x(k),x(k))$.
Then we have the following analytic expression for the image $b$:
\begin{align}\label{eq:main}
x(k) = \tau_0(k)-\tau_0(k-1)-\tau_1(k)+\tau_1(k-1).
\end{align}
Since the time evolution of the box-ball system is linearlized on the set
of the rigged configurations, this result gives an explicit solution
for the initial value problem of the box-ball systems.

\begin{proof}[Sketch of the proof of (\ref{eq:main})]
The main step of the proof is to show the following
interpretation of the tau functions.
For the given path $b=b_1\otimes b_2\otimes\cdots$, define
$T_\infty(b)=b_1^{(1)}\otimes b_2^{(1)}\otimes\cdots$,
$T_\infty^2(b)=b_1^{(2)}\otimes b_2^{(2)}\otimes\cdots$, and so on.
Then we have to show the following interpretation:
\begin{align}
\nonumber
\tau_r(k)=&\,(1-r)\times(\mbox{number of balls in }b_1\otimes b_2\otimes\cdots\otimes b_k)\\
&+\sum_{i\geq 1}
(\mbox{number of balls in }b_1^{(i)}\otimes b_2^{(i)}\otimes\cdots\otimes b_k^{(i)}).
\label{tau=rho}
\end{align}
For example, in the example of Section \ref{sec:bbs},
we have $\tau_0(8)=9$ and $\tau_1(8)=5$.
Since balls always move rightwards,
the summation in the second term is always finite.
From (\ref{tau=rho}) we can easily deduce (\ref{eq:main}).

Proof of (\ref{tau=rho}) proceeds as follows.
From the expression (\ref{def:tau}) we can construct determinants
from which we obtain the tau functions $\tau_r$ as the ultradiscrete limit.
Then by using a calculus of determinants we can show that the
tau functions satisfy the ultradiscrete Hirota bilinear form.
The Hirota bilinear form implies that the functions $\tau_r$ corresponds to
the same dynamics of the box-ball systems.
Unfortunately this is not the whole story.
The main difficulty is the fact that the analytic expression in (\ref{def:tau})
is very different from the combinatorial definition of the map $\Phi^{-1}$
and thus it is quite difficult to compare.

We do this in the following way.
The proof is induction on the rank $n$ of $A^{(1)}_n$.
Since we know that the tau functions satisfy the same dynamics of
the box-ball systems, it is enough to consider a state $T_\infty^N(b)$ where $N\gg1$.
We call such a state the asymptotic state.
Since we have the inverse scattering formalism which is the consequence
of the most part of the paper \cite{KSS}, we can easily obtain the corresponding
asymptotic rigged configuration.
Then we invoke the result of \cite{Sak1} to reduce the problem to the case of $A^{(1)}_{n-1}$
(the case $A^{(1)}_1$ can be shown by \cite{Sak1}).
This part is logically a bit complicated and we will omit the details.
Here we remark that we use the fact that the tau functions for
the general $A^{(1)}_n$ have a similar recursive structure with respect to the rank
and that we use the Yang--Baxter relations for the affine combinatorial $R$-matrices
to represent the right hand side of (\ref{tau=rho}) by the energy function and the combinatorial $R$-matrices.
Thus the proof heavily utilizes the infinite dimensional symmetry behind the box-ball system.
\end{proof}

Finally we remark that the conserved quantities $E_l$ of \cite{FOY}
indeed coincide with the rigged configurations \cite{Sak2}:
\begin{align}\label{E=Q}
E_l(b)=Q_l(\nu)
\end{align}
where $\Phi(b)=(\nu,J)$.\footnote{In \cite{Takagi2}, Takagi introduced a
scheme to factorize the dynamics of the box-ball systems of type $A^{(1)}_n$
into $A^{(1)}_1$ case by using the time evolution corresponding to the carrier of type $B^{2,1}$.
This scheme is rephrased into the rigged configuration language \cite[Section 2.7]{KOSTY}
to factorize the map $\Phi$ for general $A^{(1)}_n$ case into the map $\Phi$ for $A^{(1)}_1$ case
by using the $B^{2,1}$ type time evolution.
The proof of (\ref{E=Q}) in \cite{Sak2} uses a refinement of the latter result.}
There is a generalization of this formula for the most general rigged configuration bijection of type $A^{(1)}_n$
(see Section \ref{sec:further}).
By (\ref{E=Q}) we see that the configurations have a simple
representation theoretical origin in terms of the combinatorial $R$-matrices and
the energy function.
In fact, in \cite{Sak2} it is shown that the quantities $E_l$ (and suitable refinements)
give enough information to reconstruct riggings.
In this sense, the combinatorial algorithm of $\Phi$ itself has
a representation theoretical origin via the time evolutions of the box-ball systems.
However the mathematical origin of the riggings is still unclear.
In fact, the riggings depend on the information about the final positions where the corresponding
string is finally changed during the procedure $\Phi$.
This information is rather combinatorial and we cannot get rid of its difficulty
even if we use the representation theoretical interpretation of
the algorithm of $\Phi$ discussed above.

\section{Interlude: the box-basket-ball systems}\label{sec:bbbs}
In this section, we explain the basic properties of the {\bf box-basket-ball systems}
(BBBS for short) introduced by \cite{LPS}.
The starting point of the construction is to replace the combinatorial $R$-matrices
in the definition of the box-ball system by the whurl relations of \cite{LP}.
Rather non-trivially, the resulting dynamical system becomes a soliton system.
The characteristic property of the BBBS is that the system contains
the fermions (balls) and bosons (baskets) with mutual interaction between them.
We remark here that the BBBS is different from the super-symmetric box-ball system
of \cite{HI} constructed from the crystals for the quantum superalgebra \cite{BKK}
since their system is the extension of the box-ball system by adding another kind
of the fermionic particles.

In order to obtain intuition about the model, it is convenient to
start from a combinatorial description of the time evolution of the BBBS.
Let $b=b_1\otimes b_2\otimes\cdots$ be the state of the BBBS.
In this situation, each state is parametrized as a three dimensional vector
$(a,b,c)\in\mathbb{Z}^3$.
Our interpretation of each parameter is as follows;
$b$ is the number of {\bf baskets}, $c$ is the number of {\bf balls}
and $a$ is the number of {\bf empty places} that can fit extra balls.
The meaning of such an interpretation will become clear
when we explain the combinatorial description of the time evolution.

In the rest of this section, we consider the following situation.
We put many capacity one boxes on a line.
If necessary, we put enough many empty boxes on the right of the state.
As the rule, each box or basket can accommodate at most one ball
whereas we can put more than one baskets on a box.
Thus the balls are {\bf fermionic} particles and the baskets are {\bf bosonic} particles.
There is a nontrivial interaction between the two kinds of particles
by placing a ball within a basket.
If necessary we assume that a ball is always placed in a box
before placed in a basket.
We introduce several definitions that will be used later.
Let $V = (1,0,0)$, $F =(0,0,1)$, $B_i =(i+1,i,0)$ and $U_i =(i,i,1)$ where $i \geq 1$.
Here we give several diagrams that represent these symbols:
\begin{center}
\unitlength 20pt
\begin{picture}(14,1.7)
\put(0,0.15){$V\,\,=$}
\put(1.8,0){\hako}
\put(2.8,0){,}
\put(3.5,0.15){$F\,\,=$}
\put(5.2,0){\hako}
\put(5.2,0){\tama}
\put(6.2,0){,}
\put(7,0.15){$B_2\,\,=$}
\put(9,0){\hako}
\multiput(9,1)(0,0.5){2}{\kago}
\put(10,0){,}
\put(10.8,0.15){$U_2\,\,=$}
\put(12.8,0){\hako}
\put(12.8,0){\tama}
\multiput(12.8,1)(0,0.5){2}{\kago}
\put(13.8,0){.}
\end{picture}
\end{center}

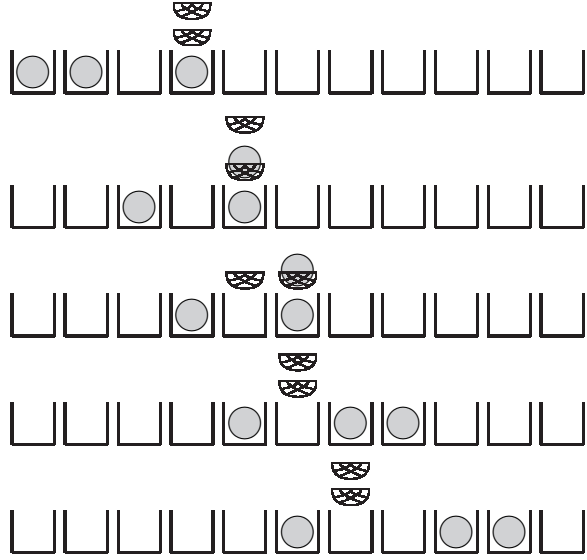
\begin{figure}
\begin{center}
\unitlength 20pt
\begin{picture}(11,2.0)
\put(0,0){\tama}
\put(1,0){\tama}
\put(3,0){\tama}
\multiput(0,0)(1,0){11}{\hako}
\multiput(3,0.9)(0,0.5){2}{\kago}
\end{picture}\\
\begin{picture}(11,2.5)
\put(2,0){\tama}
\put(4,0){\tama}
\put(4,0.85){\tama}
\multiput(0,0)(1,0){11}{\hako}
\multiput(4,0.9)(0,0.9){2}{\kago}
\end{picture}\\
\begin{picture}(11,2.0)
\put(3,0){\tama}
\put(5,0){\tama}
\put(5,0.85){\tama}
\multiput(0,0)(1,0){11}{\hako}
\put(4,0.9){\kago}
\put(5,0.9){\kago}
\end{picture}\\
\begin{picture}(11,2.0)
\put(4,0){\tama}
\put(6,0){\tama}
\put(7,0){\tama}
\multiput(0,0)(1,0){11}{\hako}
\multiput(5,0.9)(0,0.5){2}{\kago}
\end{picture}\\
\begin{picture}(11,2.0)
\put(5,0){\tama}
\put(8,0){\tama}
\put(9,0){\tama}
\multiput(0,0)(1,0){11}{\hako}
\multiput(6,0.9)(0,0.5){2}{\kago}
\end{picture}
\end{center}
\caption{Example of the time evolution of the BBBS.}
\label{fig:bbbs}
\end{figure}

Now we explain the time evolution rule.
We start from an initial state that contains finitely many baskets and balls.
\begin{quotation}
\noindent
{\bf The time evolution of the BBBS:}
First, move every empty basket to the right one step.
Full baskets are not moved.  Second, consider each ball from left to right
and move the ball to the next available empty box or basket.
Each ball is moved exactly once.
\end{quotation}
Note that if there is no basket, the above rule coincides with
the one for the box-ball system.
We give a simple but nontrivial example in Figure \ref{fig:bbbs}.

The BBBS can be constructed from the whurl relation
\begin{align}
R:(a,b,c)\otimes (d,e,f)\mapsto
(d',e',f')\otimes (a',b',c')
\end{align}
where the explicit relations are
\begin{align}
a'&= a - \min(a+b, a+c, b+f) + \min(e+c, d+c, d+b)\nonumber\\
b'&= b - \min(a+b, a+c, b+f) + \min(a+e, d+f, e+f)\nonumber\\
c'&= c - \min(e+c, d+c, d+b) + \min(a+e, d+f, e+f)\nonumber\\
d'&= d + \min(a+b, a+c, b+f) - \min(e+c, d+c, d+b)\nonumber\\
e'&= e + \min(a+b, a+c, b+f) - \min(a+e, d+f, e+f)\nonumber\\
f'&= f + \min(e+c, d+c, d+b) - \min(a+e, d+f, e+f).
\end{align}
We apply the whurl relation with $u_l:=(l,0,0)$ to the diagram
(\ref{def:T_l}) to define the operators $T_l$ ($l\in\mathbb{Z}_{l\geq 1}$).
Then we can show that the above mentioned combinatorial definition of the time evolution
coincides with $T_\infty$.
Since the whurl relations satisfy the Yang--Baxter relation,
we can show that $T_lT_k(b)=T_kT_l(b)$ for arbitrary $l,k\in\mathbb{Z}_{\geq 1}$ and states $b$.
Thus the BBBS possesses the quantum integrability.
We remark that since we do not know the underlying symmetry of the whurl relations,
we have not been able to construct a conserved quantity analogous to $E_l$.

Even if we know that the BBBS is a quantum integrable system,
it is far from clear whether the system is a soliton system.
Below we explain that the system is indeed solitonic.
For this purpose we classify solitary waves which do not change their shapes
during the free propagations under $T_\infty$.
As the result, we see that there are the following two cases.
\begin{enumerate}
\item
A consecutive sequence of $k$ balls $F_k:=FF\cdots F$.
Under the free propagation by $T_\infty$, $F_k$ moves at velocity $k$.
\item
Any sequence of $F,B,U$ which does not contain the consecutive subsequence $FF$ or $FU$,
which we call a {\bf slow soliton}.
Under the free propagation by $T_\infty$, the slow solitons move at velocity 1.
\end{enumerate}
Note that $F_k$ are the usual solitons of the box-ball system
whereas the slow solitons are the new feature of the BBBS.

Let us clarify what are the slow solitons.
The answer comes from the analysis of the {\bf phase shift}.
Here the meaning of the phase shift is as follows.
Let $A$ and $B$ are solitons on a line and suppose that they make collision
during the time evolution and retain their original form after the collision.
Then we compare the position of the soliton after the collision with
the position of the corresponding soliton supposing that there is no collision.
This difference (rightwards shift is positive) gives the phase shift.
We summarize the basic physical properties of the fermionic solitons $F_k$
and the bosonic solitons $B_{a_1}B_{a_2}\cdots B_{a_m}$ in the following table.
\begin{center}
\begin{tabular}{|c|c|c|}
\hline
&$F_k$&$B_{a_1}B_{a_2}\cdots B_{a_m}$\\
\hline\hline
velocity&$k$&1\\
\hline
phase shift&$-2k$&$-1$\\
\hline
\end{tabular}
\end{center}
Here the phase shift is defined by the scattering with $F_l$ ($l>k$).
For example, in the example in Section \ref{sec:bbs},
we see that the length one soliton $F_1$ get shifted by $-2$
after the collision with $F_3$.\footnote{Let us mention the generalizations
to the box-ball systems of the vector representations for types $A^{(1)}_n$.
Then it is known that the phase shift coincides with the energy function
(with a different normalization) between two solitons \cite{FOY}.
Here we identify freely propagating solitons with the semistandard tableaux
and regard them as the elements of crystals $B^{1,s}$ of types $A^{(1)}_{n-1}$.
Note that since we are neglecting all 1's (empty places),
we have $A^{(1)}_{n-1}$ here.
In \cite{Sak1}, it is generalized to include all cases considered in \cite{HHIKTT,FOY}
and the scatterings of soitons are identified with the affine combinatorial $R$-matrices
(\ref{eq:affine}) where each soliton corresponds to the truncated rigged configurations.}

Let us look at two solitons $F_1$ and $B_i$ of velocity one.
If we consider the scattering with $F_k$ ($k>1$), they get shifted
by $-2$ and $-1$, respectively.
To summarize, $F_1$ and $B_i$ have the same velocity whereas they have
different values of the phase shift.
Thus during the time evolutions we may have superposition of such states
and this is the origin of the slow solitons.
Therefore, in order to analyze the slow solitons, we make scatterings with many $F_k$'s
and decompose them into elementary solitons $F_1$ and $B_i$.
For example, the example in Figure \ref{fig:bbbs} shows the decomposition of
the slow soliton $U_2$ into two elementary solitons $F_1$ and $B_2$.

Based on these observations, we define the solitons of the states of the BBBS
as the elementary solitons $F_l$ and $B_i$ which we can obtain by scattering
with many additional $F_k$'s.
Let us define the amplitudes of $F_l$ and $B_i$ by $l$ and $i$, respectively.
Then we can show that the number and amplitudes of the solitons are
preserved during the time evolution of the BBBS.
Moreover we can show that scatterings of multiple solitons can be
decomposed into two body scatterings.
Hence we see that the BBBS is solitonic.

\section{Generalizations and further developments}\label{sec:further}
In the most of the present note, we only think about the simplest possible case,
namely the vector representation of type $A^{(1)}_1$.
We do so in order to provide the basic ideas without
getting into the technical complexities.
In fact, one of the nice features of our approach is its \underline{universality}.
For example, the definition of the box-ball system in (\ref{def:T_l})
has straightforward generalizations for the tensor products of
the Kirillov--Reshetikhin crystals $B^{r,s}$ for the quantum affine algebras of types other than $A^{(1)}_n$.
Here $B^{r,s}$ is the Kirillov--Reshetikhin crystals corresponding to the weight $s\Lambda_r$,
where $\Lambda_r$ is the $r$-th fundamental weight.
In this case, instead of using $u_l$ in (\ref{def:T_l}), we use the classically highest weight element of $B^{r,s}$.
Then we denote by $T^{r,s}$ the resulting time evolutions.
Again the box and ball interpretation of the time evolution
provides a nice way to get intuition about the generalized models.
For example, for the box-ball systems corresponding to the vector representations of general non-exceptional
affine algebras, there is an interpretation of $T^{1,\infty}$ in terms of particles and anti-particles
with pair creations/annihilations \cite{HKT}.
In this final section, we will give comments on the methods of the generalizations
and further properties.

\paragraph{Known extensions.}
The rigged configuration is known to have many extensions.
Indeed it is expected that such a bijection exists for the arbitrary
Kirillov--Reshetikhin crystals corresponding to general affine quantum algebras.
As mentioned in Section \ref{sec:RC}, the bijection for type $A^{(1)}_n$
is already constructed in full generalities.
Apart from this case, we have the following generalizations.
\begin{itemize}
\item
$\bigotimes B^{1,1}$ for arbitrary non-exceptional affine algebras \cite{OSS}.
\item
$\bigotimes_i B^{r_i,1}$ for type $D^{(1)}_n$ \cite{Sch1}.
\item
$\bigotimes_i B^{1,s_i}$ for type $D^{(1)}_n$ \cite{SS}.
\item
$B^{r,s}$ of type $D^{(1)}_n$ \cite{OSakaS}.
\item
$\bigotimes B^{1,1}$ for type $E^{(1)}_6$ \cite{OSano}.
\end{itemize}
We remark that the combinatorial algorithms involved in these extensions
share many common features and the philosophy which underlies these extensions is the same.
We also remark that all these results are related with
the highest weight elements of tensor products of crystals.
However, if we think about the box-ball systems we encounter the rigged configurations
for not necessarily highest weight elements.
This extension is quite natural.
Indeed the algorithm presented in Section \ref{sec:RC} does apply to
both cases without any change.
Therefore it is quite natural to consider the Kashiwara operators (analogue
of the Chevalley generators in the crystals setting) on the set of the rigged configurations.
This is achieved in \cite{Sch2} for all simply laced cases.
Remarkably, the definition of the Kashiwara operators for the all cases considered in \cite{Sch2} is uniform.

\paragraph{The method of generalizations.}
As examples of the generalizations, let us consider the cases $A^{(1)}_n$ or $D^{(1)}_n$.
Then the rigged configurations take the following form:
\begin{align}
(\nu,J)=\Bigl((\nu^{(1)},J^{(1)}), (\nu^{(2)},J^{(2)}),\cdots,(\nu^{(n)},J^{(n)})\Bigr)
\end{align}
together with the Young diagrams $\mu^{(a)}$ which is determined by the shape of
the tensor product $B=\bigotimes_iB^{r_i,s_i}$ by the following rule:
each $B^{r_i,s_i}$ in $B$ corresponds to the length $s_i$ row of $\mu^{(r_i)}$.
Note that we should consider that each $(\nu^{(a)},J^{(a)})$ corresponds to the node $a\in I_0$
of the Dynkin diagram (without the 0-node, see \cite{Kac}) for the affine algebras.
For example, the vacancy number for the present case takes the following form
\begin{align}\label{eq:def_vacancy}
P^{(a)}_\ell(\nu)
&=Q_\ell(\mu^{(a)})-2Q_\ell(\nu^{(a)})+\sum_{b\in I_0,\, b\sim a}Q_\ell(\nu^{(b)}),
\end{align}
where $a\sim b$ means that the nodes $a$ and $b$ are connected
by a single edge on the Dynkin diagram.
Note that the definition (\ref{def:P}) corresponds to the special case $a=1$ and $\mu^{(1)}=(1^k)$.
There is a nice characterization of the highest weight rigged configurations.
For the given rigged configuration $(\nu,J)$, if
\begin{align}\label{eq:characterization}
P^{(a)}_{\nu^{(a)}_i}(\nu)\geq J^{(a)}_i\geq 0
\end{align}
is satisfied by all the strings $(\nu^{(a)}_i,J^{(a)}_i)$, then $(\nu,J)$ is the highest weight rigged configuration.

The algorithm $\Phi$ is almost parallel to the definition given in Section \ref{sec:RC}
(see, for example, \cite[Appendix A]{Sak2} for details).
To get a feeling of the algorithm, suppose that we have a letter $a$ in a type $\bigotimes B^{1,1}$ path
(in Section \ref{sec:RC}, we described the case $a=2$).
Here we identify the elements of crystals $B^{r,s}$ with semistandard tableaux or Kashiwara--Nakashima
tableaux (generalizations of the semistandard tableau, see \cite{KN}).
Then we have to add a box
to each of $\nu^{(a-1)}$, $\nu^{(a-2)}$, $\cdots$, $\nu^{(1)}$ in this order.
The rule for the addition to $\nu^{(a-1)}$ is exactly the same one given in Section \ref{sec:RC}.
Suppose that we have added a box to the $\ell^{(a-1)}$-th column of $\nu^{(a-1)}$.
Then we look for the longest singular string of $(\nu^{(a-2)},J^{(a-2)})$
\underline{whose length does not exceed $\ell^{(a-1)}$} to determine where to add a box.
We do this recursively until $(\nu^{(1)},J^{(1)})$ by recursively defining $\ell^{(b)}$'s.
Finally change the riggings by using the new vacancy numbers as in Section \ref{sec:RC}.
For the modifications required for the negative letters in type $D^{(1)}_n$,
see, for example, \cite{OSakaS}.
Roughly speaking, we do an almost similar procedure twice 
(for $\bar{a}$, first proceed from $\nu^{(a)}$ to $\nu^{(n)}$
and next to the left from $\nu^{(n)}$ as above) following the crystal
graph for the vector representations $B^{1,1}$ of type $D^{(1)}_n$.

The inverse scattering formalism (\ref{eq:ist}) holds almost identically for the general cases.
For arbitrary $T^{r,s}$ of type $A^{(1)}_n$ \cite{KOSTY} and
for $T^{1,s}$ of type $D^{(1)}_n$ \cite{KSY2}, it is known that the only change caused by $T^{r,s}$
is the shift in the rigging
\begin{align}
\bigl(\nu^{(r)}_i,J^{(r)}_i\bigr)\longmapsto \bigl(\nu^{(r)}_i,J^{(r)}_i+\min(s,\nu^{(r)}_i)\bigr)
\end{align}
and all the other places do not change.
We expect that a parallel formalism should exist for all types of the quantum affine
algebras once the corresponding rigged configuration bijection is established.

Similarly, the relation (\ref{E=Q}) has the following straightforward generalization
for arbitrary rigged configurations of type $A^{(1)}_n$ (not necessarily highest weight).
In (\ref{def:T_l}) and (\ref{def:E_l}), we use the time evolution $T^{r,s}$ instead of $T_l(=T^{1,l})$.
Then we can define $E^{r,s}(b)$ as generalizations of $E_l(b)$.
Then we have \cite{Sak2}
\begin{align}
E^{r,s}(b)=Q_s(\nu^{(r)}).
\end{align}

On the other hand, the ultradiscrete tau functions formalism is only available
for the case $\bigotimes_iB^{1,s_i}$ of type $A^{(1)}_n$.
Perhaps we need to thoroughly understand the dynamics of the box-ball systems
for general cases (say, the case $\bigotimes_iB^{r_i,1}$ of type $A^{(1)}_n$).

\paragraph{Further properties.}
So far we have explained that the rigged configurations behave
very nicely with respect to the box-ball systems.
In particular, the rigged configurations have the concrete mathematical meaning
as the action and angle variables for the box-ball systems.
As the final remarks we explain there are equally remarkable properties
of the rigged configurations with respect to other mathematical problems.
In most cases, the rigged configurations behave surprisingly simply
with respect to global and deep structures of the corresponding algebras
which are usually difficult to realize.
\begin{itemize}
\item
The combinatorial $R$-matrices become trivial on the level of the rigged configurations.
If the two tensor products are isomorphic under the combinatorial $R$-matrices
$R:b\mapsto b'$, we have $\Phi(b)=\Phi(b')$.
Remind that the combinatorial $R$-matrices for general situations are
highly complicated objects.
This property is confirmed in all known cases and we expect that it is true for arbitrary
quantum affine algebras.
\item
The Sch\"{u}tzenberger involution and its generalizations become almost trivial operation
(see, for example, \cite{KSS,SS}); we take complements of all the riggings
with respect to the corresponding vacancy numbers.
Note that in this case we consider only the highest weight rigged configurations
which satisfy (\ref{eq:characterization}).
\item
In \cite{OS} a new kind of bijection $\Psi$ for the rigged configurations is introduced.
The map $\Psi$ gives one to one correspondence for the following two sets;
(i) the set of the highest weight rigged configurations for arbitrary non-exceptional quantum affine algebras
of sufficiently large rank, and (ii) the set of pairs of the highest weight rigged configurations
of type $A^{(1)}_n$ and the Littlewood--Richardson tableaux.
In experiments, we can see that the map $\Psi$ coincides with (and generalizes)
the global involution exchanging the nodes 0 and $n$ of the Dynkin diagram of the algebra
(if such an involution exists, see \cite{LOS}).
Remarkably, the construction of the algorithm $\Psi$ is quite simple and does not
depend on the choices of the corresponding non-exceptional algebras.
Indeed, the algorithm coincides with the type $A^{(1)}_n$ rigged configuration bijection
if we change left and right in the definition (as if we are using mirrors).
The Littlewood--Richardson tableaux naturally appear as the recording tableaux.
We remark that such a correspondence is very difficult to construct if we do not
use the rigged configurations (see \cite{Shi}).
\item
In \cite{Sch3} the affine Kashiwara operators for type $D^{(1)}_n$ are realized via the Dynkin
involution exchanging the nodes 0 and 1.
The realization relies on a rather nontrivial bijection between the Kashiwara--Nakashima
tableaux and a combinatorial objects called the plus-minus diagrams.
Then the involution is realized as changing columns of the plus-minus diagrams.
In \cite{OSakaS} we see that the plus-minus diagrams essentially coincide with the rigged configurations.
Thus we can realize the Dynkin involution $0\leftrightarrow 1$
as a transformation on the rigged configurations.

However the main point of the result is not the practical values.
Rather, the result reveals that the crystal structure of the corresponding case is
essentially governed by the rigged configurations.
Note that in this case the rigged configurations for non-highest weight
elements play the role.
\end{itemize}

\paragraph{Concluding words.}
We have seen that the rigged configurations have very
special properties which are usually difficult to see so that it is tempting to say that
they are one of the canonical realizations of the Kirillov--Reshetikhin crystals.
Not only they give a nice presentation, they also have concrete mathematical meanings
and it seems that they originate from deep aspects of the infinite dimensional
symmetry of the quantum affine algebras.
Although the theory of the rigged configurations is still in a very early stage,
we expect that the progress of the theory will give unique insights into
the nature of the symmetry of the quantum affine algebras.

\bigskip

\begin{flushright}
\shadowbox{Date: 12/12/12}
\end{flushright}

\end{document}